\documentclass[a4paper,9pt]{article}

\usepackage{amsmath,amscd,amssymb}

\def\l{\lambda}

\def\h{\mathfrak{h}}

\def\C{\mathbb{C}}
\def\Z{\mathbb{Z}}
\def\N{\mathbb{N}}
\def\R{\mathbb{R}}
\def\Q{\mathbb{Q}}
\def\CB{{\mathcal{C}}_B}
\def\DB{{\mathcal{D}}_B}
\def\CF{{\mathcal{C}}_F}

\newtheorem{theo}{\bf{Theorem}}
\newtheorem{prop}{Proposition}
\newtheorem{lem}{Lemma}
\newtheorem{cor}{Corollary}
\begin{document}
\title{A link between two elliptic quantum groups}
\author{Pavel Etingof and Olivier Schiffmann}
\maketitle
\begin{abstract} We consider the category $\CB$ of meromorphic finite-dimensional representations of the quantum elliptic algebra $\mathcal{B}$ constructed via Belavin's R-matrix, and the category $\CF$ of meromorphic finite-dimensional representations of Felder's elliptic quantum group $\mathcal{E}_{\tau, \frac{\gamma}{2}}(\mathfrak{gl}_n)$. For any fixed $c \in \C$, we use a version of the Vertex-IRF correspondence to construct two families of (generically) fully faithful functors $\mathcal{H}^c_{x}: \CB \to \DB$ and $\mathcal{F}^c_{x} :\CF \to \DB$ where $\DB$ is a certain category of infinite-dimensional representations of $\mathcal{B}$ by difference operators. We use this to construct an equivalence between the abelian subcategory of $\CB$ generated by tensor products of vector representations and the abelian subcategory of $\CF$ generated by tensor products of vector representations. 
\end{abstract}
\section{Categories of meromorphic representations}
\paragraph{}In this section, we recall the definitions of various categories of representations of quantum elliptic algebras.
\paragraph{Notations:} let us fix $\tau \in \C$, $\mathrm{Im}(\tau) > 0$, $\gamma \in \R\backslash \Q$ and $n \geq 2$. Denote by $(v_i)_{i=1}^n$ the canonical basis of $\C^n$ and by $(E_{ij})_{i,j=1}^n$ the canonical basis of $\mathrm{End}(\C^n)$, i.e $E_{ij}v_k=\delta_{jk}v_i$ . Let $\h=\{\sum_i \lambda_i E_{ii}\;|\; \sum_i \lambda_i=0\}$ be the space of diagonal traceless matrices. We have a natural identification $\h^*=\{\sum_i \lambda_i E_{ii}^*\;|\sum_{I}\lambda_i=0\}$. In particular, the weight of $v_i$ is $\omega_i= E_{ii}^*-\frac{1}{n}\sum_{k}E_{kk}^*$.

\paragraph{Classical theta functions:} the theta function $\theta_{\kappa,\kappa'}(t;\tau)$  with characteristics $\kappa, \kappa'\in \R$ is defined by the formula
 $$\theta_{\kappa,\kappa'}(t;\tau)=\sum_{m \in \Z} e^{i \pi (m + \kappa)((m+\kappa)\tau + 2 (t + \kappa'))}.$$
It is an entire function whose zeros are simple and form the (shifted) lattice $\{\frac{1}{2}-\kappa + (\frac{1}{2} - \kappa')\tau \} + \Z+\tau\Z$.\\
\hbox to1em{\hfill}Theta functions satisfy (and are characterized up to renormalization by) the following fundamental monodromy relations
\begin{align}
\theta_{\kappa,\kappa'}(t+1;\tau)&=e^{2i\pi\kappa}\theta_{\kappa,\kappa'}(t;\tau) \label{E:theta1},\\
\theta_{\kappa,\kappa'}(t+\tau;\tau)&=e^{-i\pi\tau-2i\pi(t+\kappa')}\theta_{\kappa,\kappa'}(t;\tau).\label{E:theta2}
\end{align}
Theta functions with different characteristics are related to each other by shifts of $t$:
\begin{equation}\label{E:theta3}
\theta_{\kappa_1+\kappa_2,\kappa'_1+\kappa'_2}(t;\tau)=e^{i\pi\kappa_2^2\tau + 2i\pi \kappa_2 (t+\kappa'_1 + \kappa'_2)}\theta_{\kappa_1,\kappa'_1}(t+\kappa_2\tau+\kappa'_2;\tau).
\end{equation}
In particular, we set $\theta (t)=\theta_{\frac{1}{2},\frac{1}{2}}(t;\tau)$.

\subsection{Meromorphic representations of the Belavin quantum elliptic algebra}
Consider the two $n \times n$ matrices
 $$A=\pmatrix 1&0&\ldots&0\\ 0&\xi&\ldots&0\\\vdots& \vdots&
\ddots&\vdots\\ 0&0&\ldots&\xi^{n-1} \endpmatrix \qquad B= \pmatrix
0&1&\ldots&0\\ \vdots& \vdots& \ddots& \vdots  \\0&0&\vdots&1\\ 1&0&\ldots &0
\endpmatrix$$
where $\xi=e^{2i\pi/n}$. We have $A^n=B^n=Id, \; BA= \xi AB$, i.e $A,B$ generate the Heisenberg group. Belavin (\cite{Bel}) introduced the matrix $R^{B}(z) \in \mathrm{End}(\C^n) \otimes \mathrm{End}(\C^n)$, uniquely determined by the following properties:
\begin{enumerate}
\item Unitarity: $R^B(z)R^B_{21}(-z)=1$,
\item $R^B(z)$ is meromorphic, with simple poles at $z=\gamma + \Z + \tau \Z$,
\item $R^B(0)=P: x \otimes y \mapsto y \otimes x$ for $x,y \in \C^n$ (permutation),
\item Lattice translation properties:
\begin{align*}
R^B(z+1)&= A_{1}R^B(z)A_{1}^{-1}=A_{2}^{-1}R^B(z)A_{2},\\
R^B(z+\tau)&= e^{-2i \pi \frac{n-1}{n} \gamma}B_{1}R^B(z)B_{1}^{-1}=e^{-2i \pi
\frac{n-1}{n} \gamma} B_{2}^{-1}R^B(z)B_{2}.
\end{align*}
\end{enumerate}
In particular, $R^B(z)$ commutes with $A \otimes A$ and $B \otimes B$. The matrix $R^B(z)$ satisfies the quantum Yang-Baxter equation with spectral parameters:
\begin{equation*}
R^B_{12}(z-w)R^B_{13}(z)R^B_{23}(w)=R^B_{23}(w)R^B_{13}(z)R^B_{12}(z-w).
\end{equation*}
 \paragraph{The category $\CB$:} following Faddeev, Reshetikhin, Takhtajan and Semenov-Tian-Shansky, one can define an algebra $\mathcal{B}$ from $R^B(z)$, using the RLL formalism-see \cite{FRT}, \cite{RS}. However, we will only need to consider a certain category of modules over this algebra, defined as follows.\\
\hbox to1em{\hfill}Let $\CB$ be the category whose objects are pairs $(V,L(z))$ where $V$ is a finite dimensional vector space and $L(z) \in \mathrm{End}(\C^n) \otimes \mathrm{End}(V)$ is an invertible meromorphic function (the L-operator) such that $L(z+n)=L(z)$ and $L(z+n \tau)=L(z)$, satisfying the following relation in the space $\mathrm{End}(\C^n) \otimes \mathrm{End}(V) \otimes \mathrm{End}(V)$:
\begin{equation}\label{E:01}
R^B_{12}(z-w)L_{13}(z)L_{23}(w)=L_{23}(w)L_{13}(z)R^B_{12}(z-w)
\end{equation}
(as meromorphic functions of $z$ and $w$); morphisms $(V,L(z)) \to (V', L'(z))$ are linear maps $\varphi: V \to V'$ such that $(1 \otimes \varphi)L(z)=L'(z)(1 \otimes \varphi)$ in the space $\mathrm{Hom}(\C^n \otimes V, \C^n \otimes V')$.  The quantum Yang-Baxter relation for $R^B$ implies that $(\C^n, \chi(z)R^B(z-w)) \in {\mathcal{O}}b(\CB)$ for all $w \in \C$, where we set $\chi(z)=\frac{\theta(z-(1-\frac{1}{n})\gamma)}{\theta(z)}$. This object is called the vector representation and will be denoted simply by $V_B(w)$.\\
\hbox to1em{\hfill}The category $\CB$ is naturally a tensor category with tensor product
\begin{equation}\label{E:tensB}
(V,L(z)) \otimes (V',L'(z))=(V \otimes V', L_{12}(z)L'_{13}(z))
\end{equation}
at the level of objects and with the usual tensor product at the level of morphisms.\\
\hbox to1em{\hfill}There is a notion of a dual representation in the category $\CB$: the (right) dual of $(V,L(z))$ is $(V^*,L^*(z))$ where $L^*(z)=L^{-1}(z)^{t_2}$ (first apply inversion, then apply the transposition in the second component $t_2$). If $V,W \in \mathcal{O}b(\CB)$ and $\varphi \in \mathrm{Hom}_{\CB}(V,W)$ then $\varphi^t \in \mathrm{Hom}_{\CB}(W^*,V^*)$.

\paragraph{}We will also need an extended category $\CB^{x}$ defined as follows: objects of $\CB^{x}$ are objects of $\CB$ but we set
$$\mathrm{Hom}_{\CB^{x}}(V,V')=\mathrm{Hom}_{\CB}(V,V') \otimes M_{\C}$$
where $M_{\C}$ is the field of meromorphic functions of a complex variable $x$. In other words, morphisms in $\CB^{x}$ are meromorphic 1-parameter families of morphisms in $\CB$.
\paragraph{The category $\DB$:} We now define a difference-operator variant of the categories $\mathcal{C}_B, \CB^{x}$. Let us denote by $M_{\h^*}$ the field of $(n\omega_i)$-periodic meromorphic functions $\h^* \to \C$ and by $D_{\h^*}$ the $\C$-algebra generated by $M_{\h^*}$ and shift operators $T_\mu:\,M_{\h^*} \to M_{\h^*}, f(\lambda) \mapsto f(\lambda+\mu)$ for $\mu \in \h^*$. If $V$ is a finite-dimensional vector space, we set $V_{\h^*}= M_{\h^*} \otimes_{\C} V$, and $D(V)= D_{\h^*} \otimes \mathrm{End}(V)$. Let $\mathcal{D}_B$ be the category whose objects are pairs $(V,L(z))$ where $V$ is a finite-dimensional $\C$-vector space and $L(z) \in \mathrm{End}(\C^n) \otimes D(V)$ is an invertible operator with meromorphic coefficients satisfying (\ref{E:01}) in $\mathrm{End}(\C^n) \otimes D(V) \otimes D(V)$; morphisms $(V, L(z)) \to (V', L'(z))$ are $(n\omega_i)$-periodic meromorphic functions $ \varphi:\;\h^* \to \mathrm{Hom} (V, V')$ such that $(1 \otimes \varphi)L(z)=L(z)(1 \otimes \varphi)$ in $\mathrm{Hom}_\C( \C^n \otimes V_{\h^*}, \C^n \otimes V'_{\h^*})$ (i.e morphisms are $M_{\h^*}$-linear).\\
\hbox to1em{\hfill}The category $\DB$ is a right-module category over $\CB$, i.e we have a (bi)functor $\otimes:\,\DB \times \CB \to \DB$ defined by (\ref{E:tensB}).
\paragraph{}The category $\DB^{x}$ is defined in an analogous way: objects are pairs $(V,L(z,x))$ as in $\DB$ but the L-operator is now a meromorphic function of $z$ and $x$, and morphisms $(V,L(z,x)) \to (V',L'(z,x))$ are meromorphic maps $\varphi(\lambda,x):\h^* \times \C \to \mathrm{Hom}_\C(V,V')$ satisfying $(1 \otimes \varphi)L(z,x)=L(z,x)(1 \otimes \varphi)$.
\subsection{Meromorphic representations of the elliptic quantum group $\mathcal{E}_{\tau,\gamma/2}(\mathfrak{gl}_n)$}
\paragraph{Felder's dynamical R-matrix:} let us consider the functions of two complex variables
 $$\alpha (z,l)= \frac{\theta (l+ \gamma)\theta(z)}{\theta(l)\theta(z-\gamma)}, \qquad \beta(z,l)=\frac{\theta(z-l)\theta(\gamma)}{\theta(l)\theta(z-\gamma)}.$$
 As functions of $z$, $\alpha$ and $\beta$ have simple poles at $z= \gamma + \Z + \tau \Z$ and satisfy
\begin{alignat*}{2}
\alpha(z+1,l) & = \alpha(z,l), & \qquad
\alpha(z+\tau,l) & = e^{-2i\pi \gamma} \alpha(z,l),  \\
 \beta(z+1,l)& = \beta(z,l), & \qquad
\beta(z+\tau,l) & = e^{-2i\pi (\gamma - l)} \beta(z,l).
\end{alignat*}
\hbox to1em{\hfill}Felder introduced in \cite{F} the matrix $R^{F}(z,\l): \C \times \h^* \to \mathrm{End}(\C^{n}) \otimes \mathrm{End}(\C^{n})$:
$$R^{F}(z,\l)= \sum_{i}E_{ii}\otimes E_{ii} + \sum_{i \neq j} \alpha(z,
\l_{i}-\l_{j})E_{ii}\otimes E_{jj} +
\sum_{i \neq j} \beta(z, \l_{i}-\l_{j})E_{ji} \otimes E_{ij}$$
where $\l=\sum_i\l_{i}E_{ii}^* \in \h^*$.\\
\hbox to1em{\hfill}This matrix is a solution of the quantum dynamical Yang-Baxter equation with spectral parameters
\begin{align*}
R^F_{12}(z-w, \l -\gamma h_{3})&R^F_{13}(z, \l)R^F_{23}(w,\l-\gamma h_1)\\
=&R^F_{23}(w,\l)R^F_{13}(z, \l -\gamma h_2)R^F_{12}(z-w,\l)
\end{align*}
where we have used the following convention: if $V_i$ are diagonalizable $\h$-modules with weight decomposition $V_i=\bigoplus_\mu V_i^{\mu}$ and $a(\lambda) \in \mathrm{End}(\bigotimes_iV_i)$ then
$$a(\lambda - \gamma h_l)_{|\bigotimes_i V_i^{\mu_i}}=a(\lambda-\gamma \mu_l)$$
 As usual, indices indicate the components of the tensor product on which the operators act.\\
\hbox to1em{\hfill}In addition, $R(z,\lambda)$ satisfies the following two conditions:
\begin{enumerate}
\item Unitarity: $R_{12}(z, \l)R_{21}(-z,\l)=Id,$
\item Weight zero: $ \forall h \in \h,\;[h^{(1)}+h^{(2)}, R(z,\l)]=0.$
\end{enumerate}
\paragraph{The category $\CF$:} It is possible to use $R(z,\lambda)$ to define an algebra by the RLL-formalism (see \cite{F}): the elliptic quantum group $\mathcal{E}_{\tau, \gamma/2}(\mathfrak{gl}_n(\C))$. However, we will only need the following category of its representations $\CF$, introduced by Felder in \cite{F} and studied by Felder and Varchenko in \cite{FV1}: objects are pairs $(V,L(z,\lambda))$ where $V$ is a finite-dimensional diagonalizable $\h$-module and $L(z,\lambda): \C \times \h^* \to \mathrm{End}(\C^n) \otimes \mathrm{End}(V)$ is an invertible meromorphic function which is $(n\omega_i)$-periodic in $\lambda$ and which satisfies the following two conditions:
$$ [h_1+h_2,L(z,\l)]=0,$$ 
\begin{align}\label{E:02}
R_{12}(z -w, \l -\gamma h_3)&L_{13}(z,\l)L_{23}(w, \l-\gamma h_1)
\notag \\
=&L_{23}(w, \l)L_{13}(z, \l -\gamma h_2)R_{12}(z-w, \l)
\end{align}
Morphisms $(V,L(z,\lambda)) \to (V',L'(z, \lambda))$ are $(n\omega_i)$-periodic meromorphic weight zero maps $\varphi(\lambda):V \to V'$ such that $L'(z,\lambda) ( 1 \otimes \varphi(\lambda -\gamma h_1))=(1 \otimes \varphi(\lambda))L(z,\lambda)$. The dynamical quantum Yang-Baxter relation for $R^F(z,\lambda)$ implies that $(\C^n,\\ R^F(z-w,\lambda)) \in {{\mathcal{O}}b}(\CF)$ for all $w \in \C$. This is the vector representation and it will be denoted by $V_F(w)$.\\
\hbox to1em{\hfill}The category $\CF$ is naturally equipped with a tensor structure: it is defined on objects by 
$$(V,L(z, \lambda)) \otimes (V',L'(z,\lambda))=(V \otimes V', L_{12}(z,\lambda-\gamma h_3)L'_{13}(z,\lambda)),$$ 
and if $\varphi \in \mathrm{Hom}_{\CF}(V,W), \varphi' \in \mathrm{Hom}_{\CF}(V',W')$ then
$$(\varphi \otimes \varphi')(\lambda)=\varphi(\lambda - \gamma h_2) \otimes \varphi'(\lambda) \in \mathrm{Hom}_{\CF}(V \otimes V',W \otimes W').$$\\
\hbox to1em{\hfill}There is a notion of a dual representation in the category $\CF$: the (right) dual of $(V,L(z,\lambda))$ is $(V^*,L^*(z,\lambda))$ where $L^*(z,\lambda)=L^{-1}(z,\lambda +\gamma h_2)^{t_2}$ (apply inversion, shifting and then apply the transposition in the second component $t_2$). If $V,W \in \mathcal{O}b(\CB)$ and $\varphi(\lambda) \in \mathrm{Hom}_{\CB}(V,W)$ then $\varphi^*(\lambda):=\varphi(\lambda + \gamma h_1)^t \in \mathrm{Hom}_{\CB}(W^*,V^*)$.
\paragraph{}The extended category $\CF^{x}$ is defined by $\mathcal{O}b(\CF^{x})=\mathcal{O}b(\CF)$ and
$$\mathrm{Hom}_{\CF^{x}}(V,V')=\mathrm{Hom}_{\CF}(V,V')\otimes M_{\C}$$
i.e morphisms in $\CF^{x}$ are meromorphic 1-parameter families of morphisms in $\CF$.
\section{The functor $\mathcal{F}^c_x:\CF \to \DB$}
\paragraph{}In this section, we define a family of functors from meromorphic (finite-dimensional) representations of $\mathcal{E}_{\tau, \frac{\gamma}{2}}(\mathfrak{gl}_n(\C))$ to infinite-dimensional representations of the quantum elliptic algebra $\mathcal{B}$.
\subsection{Twists by difference operators:}
 \paragraph{}For any finite-dimensional diagonalizable $\h$-module $V$, let $e^{\gamma D}\in \mathrm{End}(V_{\h^*})$ denote the shift operator: $e^{\gamma D}. \sum_\mu f_\mu (\lambda) v_\mu=\sum_\mu f(\lambda + \gamma \mu) v_\mu$, $v_\mu \in V_\mu$. Now let $(V, L(z, \lambda)) \in \CF$, and let $S(z, \lambda), S'(z,\lambda): \C \times \h^* \to \mathrm{End}(\C^n)$ be meromorphic and nondegenerate. Define the difference-twist of $(V,L(z,\lambda))$ to be the pair $(V,L^{S,S'}(z))$ where
\begin{equation}\label{E:03}
L^{S,S'}(z)=S_1(z,\lambda-\gamma h_2) L(z,\lambda)e^{-\gamma D_1} S'_1(z,\lambda)^{-1} \in \mathrm{End}(\C^n) \otimes D(V).
\end{equation}
This is a difference operator acting on $\C^n \otimes V_{\h^*}$.
\begin{lem} The difference operator $L^S(z,\lambda)$ satisfies the following relation in $\mathrm{End}(\C^n) \otimes D(V) \otimes D(V)$:
$$T_{12}(z,w,\lambda-\gamma h_3) L_{13}^{S,S'}(z)L_{23}^{S,S'}(w)=L^{S,S'}_{23}(w)L^{S,S'}_{13}(z)T_{12}'(z,w,\lambda)$$
where 
\begin{align}
T(z,w,\lambda)&=S_{2}(w,\lambda)S_1(z,\lambda-\gamma h_2)R_{12}^F(z-w,\lambda)S_2(w,\lambda-\gamma h_1)^{-1}S_1(z,\lambda)^{-1} \label{E:04}\\
T'(z,w,\lambda)&=S'_{1}(z,\lambda)S'_2(w,\lambda+\gamma h_1)R_{12}^F(z-w,\lambda)S'_1(z,\lambda+\gamma h_1)^{-1}S'_2(w,\lambda)^{-1} \label{E:05}
\end{align}
\end{lem}
\textbf{Proof:} the proof is straightforward, using relation (\ref{E:02}) for $L(z,\lambda)$ and the weight zero property of $R^F(u,\lambda)$ and $L(u,\lambda)$.$\square$
\subsection{The Vertex-IRF transform}
\paragraph{}Let $\phi_l(u)=e^{2i\pi(\frac{l^2\tau}{n}+\frac{lu}{n})}\theta_{0,0}(u+l\tau;n\tau)$ for $l=1,\ldots n$. Then the vector $\Phi(u)=(\phi_1(u), \ldots,\phi_n(u))$ is, up to renormalization, the unique holomorphic vector in $\C^n$ satisfying the following monodromy relations:
\begin{align}
\Phi(u+1)&=A\Phi(u),\label{E:06}\\
\Phi(u+\tau)&=e^{-i\pi \frac{\tau}{n}-2i\pi\frac{u}{n}}B \Phi(u) \label{E:07}
\end{align}
Now let $S(z,\lambda): \C \times \h^* \to \mathrm{End}(\C^n)$ be the matrix whose columns are $(\Phi_1(z,\lambda),\\ \ldots ,\Phi_n(z,\lambda))$ where $\Phi_j(z,\lambda)=\Phi(z-n\lambda_j)$. Using (\ref{E:06})-(\ref{E:07}), it is easy to see that we have $\mathrm{det}(S(z,\lambda))=\mathrm{Const}(\lambda)\theta(z)$ and hence that $S(z,\lambda)$ is invertible for $z \neq 0$ and generic $\lambda$.
\begin{lem} We have
\begin{align*}
R^B(z-w)S_1(z,\lambda)S_2(w,\lambda-\gamma h_1)&=S_2(w,\lambda)S_1(z,\lambda-\gamma h_2)R^F(z-w,\lambda)\\
R^B(z-w)S_2(w,\lambda)S_1(z,\lambda+\gamma h_2)&=S_1(z,\lambda)S_2(w,\lambda+\gamma h_1)R^F(z-w,\lambda)
\end{align*}
\end{lem}
\textbf{Proof:} the first relation is equivalent to the following identities for $i,j=1,\ldots n$:
\begin{align*}
R^B(z-w) \Phi_i(z,\lambda)\otimes \Phi_i(w,\lambda-\gamma \omega_i)&=\Phi_i(z,\lambda-\gamma \omega_i) \otimes \Phi_i(w,\lambda)\\
R^B(z-w)\Phi_i(z,\lambda) \otimes \Phi_j(w,\lambda-\gamma \omega_i)&= \alpha(z-w,\lambda_i-\lambda_j)\Phi_i(z,\lambda-\gamma\omega_j)\otimes \Phi_j(w,\lambda)\\
&\; + \beta(z-w,\lambda_i-\lambda_j)\Phi_j(z,\lambda-\gamma\omega_i)\otimes \Phi_i(w,\lambda)
\end{align*}
 These identities are proved by comparing poles and transformation properties under lattice translations as functions of $z$ and $w$, and using the uniqueness of $\Phi$. The second relation of the lemma is proved in the same way. These identities are essentially the Vertex/Interaction-Round-a-Face transform of statistical mechanics (see \cite{H2},\cite{KrZ} and \cite{FV2} for the case $n=2$).$\square$
\subsection{Construction of the functor $\mathcal{F}^c_x:\mathcal{C}_F \to \mathcal{D}_B$}
\paragraph{}Let us fix some $c \in \C$. We can now define the family of functors $\mathcal{F}^c_{x}:\,\CF \to \CB$ indexed by $x \in \C$: for $(V,L(z,\lambda)) \in \CF$, set $\mathcal{F}^c_{x}((V,L(z,\lambda)))=(V,L^{S_x,S_{x+c}}(z))$ with $S_u(z,\lambda)=S(z-u,\lambda)$ as above and let $\mathcal{F}^c_{x}$ be trivial at the level of morphisms.
\begin{prop} $\mathcal{F}^c_{x}:\CF \to \DB$ is a functor.
\end{prop}
\textbf{Proof:} it follows from Lemma 2 that $(V,L^{S_x,S_{x+c}}(z)) \in {\mathcal{O}}b(\DB)$. Furthermore, if $\varphi(\lambda) \in \mathrm{Hom}_{\mathcal{C}_F} ((V,L(z,\lambda)),(V',L'(z,\lambda)))$ then by definition we have $L'(z,\lambda)\\(1 \otimes \varphi(\lambda-\gamma h_1))=(1 \otimes \varphi(\lambda))L(z,\lambda)$, so that 
\begin{align*}
S_1(z-x,&\lambda-\gamma h_2)L'(z,\lambda)e^{-\gamma D_1} S_1(z-x-c,\lambda)^{-1} (1 \otimes \varphi(\lambda))\\&= S_1(z-x,\lambda - \gamma h_2)L'(z,\lambda)(1 \otimes \varphi(\lambda-\gamma h_1))e^{-\gamma D_1} S_1(z-x-c,\lambda)^{-1}\\&=S_1(z-x,\lambda - \gamma h_2)(1 \otimes \varphi(\lambda)) L'(z,\lambda) e^{-\gamma D_1} S_1(z-x-c,\lambda)^{-1}\\&= (1 \otimes \varphi(\lambda)) S_1(z-x,\lambda-\gamma h_2)L'(z,\lambda)e^{-\gamma D_1}S_1(z-x-c,\lambda)^{-1}
\end{align*}
since $\varphi(\lambda)$ is of weight zero. Thus $\mathcal{F}^c_{x}(\varphi(\lambda))$ is an intertwiner in the category $\DB$.$\square$
\paragraph{}We can also think of the family of functors $\mathcal{F}^c_{x}$ as a single functor $\mathcal{F}^c:\CF^{x} \to \DB^{x}$.
\paragraph{Remark:} we can think of the difference-twist and the relations in Lemma 2 as a dynamical analogue of the notion of equivalence of R-matrices due to Drinfeld and Belavin-see \cite{BD}. 
\section{The image of the trivial representation and the functor $\mathcal{H}^c_x: \CB \to \DB$}
\paragraph{}Applying the functor $\mathcal{F}^c_{x}$ to the trivial representation $(\C,\mathrm{Id}) \in {\mathcal{O}}b(\CF)$ yields
$$\mathcal{F}^c_{x}((\C,\mathrm{Id}))=(\C,S(z-x,\lambda)e^{-\gamma D_1} S(z-x-c,\lambda)^{-1}).$$
We will denote this object by $I^c_x$. For instance, when $n=2$, we obtain a representation of the Belavin quantum elliptic algebra as difference operators acting on the space of periodic meromorphic functions in one variable $\lambda$, i.e given by an L-operator
$$L(z)=\pmatrix a(z) & b(z)\\ c(z) & d(z) \endpmatrix$$
where $a(z),b(z),c(z),d(z)$ are operators of the form $f(z)T_{-\gamma}+g(z)$ where $T_{-\gamma}$ is the shift by $-\gamma$.
\paragraph{}Such representations of $\mathcal{B}$ by difference operators already appeared in the work of Krichever, Zabrodin (\cite{KrZ}) (for $n=2$) and Hasegawa (\cite{H1},\cite{H2})(for the general case), where they were also derived by some Vertex-IRF correspondence.
\paragraph{Definition:} Let $c \in \C$ and let $\mathcal{H}^c_{x}: \CB \to \DB$ be the functor defined by the assignment $V \to I^c_{x} \otimes V$ and which is trivial at the level of morphisms. The family of functors $\mathcal{H}^c_{x}$ gives rise to a functor $\mathcal{H}^c:\CB^{x} \to \DB^{x}$.
\section{Full Faithfulness of the functor $\mathcal{H}^c_{x}:\CB \to \DB$}
In this section, we prove the following result
\begin{prop}\label{P:02}Let $V,V' \in \mathcal{O}b(\CB)$. Then for all but finitely many values of $x \;\mathrm{mod}\; \Z + \Z\tau$, the map 
$$\mathcal{H}^c_{x}: \mathrm{Hom}_{\CB}(V,V') \stackrel{\sim}{\to} \mathrm{Hom}_{\DB}(\mathcal{H}^c_{x}(V), \mathcal{H}^c_{x}(V'))$$
is an isomorphism.
\end{prop}
\paragraph{Proof:} since $\mathrm{Hom}_{\CB}(V,V') \simeq \mathrm{Hom}_{\CB} (\C, V' \otimes V^*)$, $\mathrm{Hom}_{\DB}(I^c_{x} \otimes V, I^c_{x} \otimes V') \simeq \mathrm{Hom}_{\DB}(I^c_{x}, I^c_{x} \otimes V' \otimes V^*)$, it is enough to show that the map $\mathcal{H}^c_{x}: \mathrm{Hom}_{\CB} (\C, W) \to \mathrm{Hom}_{\DB}(I^c_{x}, I^c_{x} \otimes W)$ is an isomorphism for all $W \in \mathrm{O}b(\CB)$. Since $\mathcal{H}^c_{x}$ is trivial at the level of morphisms, this map is injective. Now let $W \in \mathcal{O}b(\CB)$ and let $\varphi(\lambda) \in \mathrm{Hom}_{\DB}(I^c_{x} , I^c_{x} \otimes W)$, that is, $\varphi(\lambda)$ is a $(n \omega_i)$-periodic meromorphic function $\h^* \to W$ satisfying the equation
\begin{equation*}
\begin{split}
\varphi_2(\lambda) S_1(z-x,\lambda)&e^{-\gamma D_1} S_1(z-x-c,\lambda)^{-1}\\
&=S_1(z-x,\lambda)e^{-\gamma D_1}S_1(z-x-c,\lambda)^{-1}L_{12}(z)\varphi_2(\lambda)
\end{split}
\end{equation*}
where $L(z)$ is the L-operator of $W$. This is equivalent to
\begin{equation}\label{E:09}
L_{12}(z)\varphi_2(\lambda)=S_1(z-x-c,\lambda)\varphi_2(\lambda + \gamma h_1)S_1(z-x-c,\lambda)^{-1}
\end{equation}
\hbox to1em{\hfill}Now $L(z)$ is an elliptic function (of periods $n$ and $n\tau$) so it is either constant or it has a pole. Restricting $W$ to the subrepresentation $\mathrm{Span}(\varphi(\lambda), \lambda \in \h^*)$, we see that the latter case is impossible for generic $x$ as the RHS of (\ref{E:09}) has a pole at $z=x+c$ only; hence $L(z)$ is constant. Furthermore, from (\ref{E:09}) we see that the matrix 
$$M(\lambda)=S_1(z-x-c,\lambda)^{-1}L_{12}S_1(z-x-c,\lambda)$$
 is independent of $z$. In particular, setting $z \mapsto z+1$ and using the transformation properties (\ref{E:06}) of $S(z,\lambda)$, we obtain $[A_1,L_{12}]=0$. This implies that $L=\sum_i E_{ii}\otimes D_i$ for some $D_i \in \mathrm{End}(W)$.
\begin{lem} Let $U$ be a finite dimensional vector space, let $T\in \mathrm{End}(\C^n) \otimes \mathrm{End}(U)$ be an invertible solution of the equation
\begin{equation}\label{E:QYB}
R^B_{12}(z)T_{13}T_{23}=T_{23}T_{13}R^B_{12}(z)
\end{equation}
such that $T=\sum_i E_{ii} \otimes D_i$ for some $D_i \in \mathrm{End}(U)$. Then $[D_i,D_j]=0$ for all $i,j$ and there exists $X \in \mathrm{End}(U)$ such that $X^n=1$ and $D_{i+1}=XD_i$ for all $i=1, \ldots n$.
\end{lem}
\textbf{Proof:} let us write $R^B(z)=\sum_{p,q,r,s}R_{p,q,r,s}(z)E_{pq}\otimes E_{rs}$. Then equation (\ref{E:QYB}) is equivalent to $R_{p,q,r,s}(z)D_pD_q=R_{p,q,r,s}(z)D_sD_r$ for all $p,q,r,s$. But it follows from the general formula for $R^B(z)$ that $R_{p,q,r,s}(z) \neq 0$ if and only if $p+q \equiv r+s \;(\mathrm{mod}\;n)$. Thus we have $[D_i,D_j]=0$ for all $i,j$ and  $X:=D_iD_{i+1}^{-1}$ is independent of $i$, and satisfies $X^n=1$.$\square$\\
\hbox to1em{\hfill}By the above lemma, there exists $X \in \mathrm{End}(W)$ such that $X^n=1$ and $D_{i+1}=XD_i$. Suppose that $X \neq 1$ and choose $e \in W$ such that $X(e)=\xi^ke$ with $\xi^k \neq 1$. Now we apply the transformation $z \mapsto z+\tau$ to the matrix $M(\lambda)$. Noting that, by (\ref{E:07}), $S(z-x-c+\tau,\lambda)=e^{-i\pi\tau/2 -2i\pi(z-x-c)/n}B S(z-x-c,\lambda)F(\lambda)$ where $F(\lambda)=\mathrm{diag}(e^{-2i\pi \lambda}, \ldots e^{-2i\pi\lambda_n})$, we obtain the equality
\begin{equation*}
\begin{split}
F(\lambda)^{-1}S_1(z-x-c,\lambda)^{-1}B_1^{-1}L_{12}B_1&S_1(z-x-c,\lambda)F(\lambda)\\
&=S_1(z-x-c,\lambda)^{-1}L_{12}S_1(z-x-c,\lambda)
\end{split}
\end{equation*}
Applying this to the vector $e$ yields $\mathrm{Ad}F(\lambda)(M(\lambda))(e)=\xi^{-k} M(\lambda)(e)$. This is possible for all $\lambda$ only if $k\equiv 0\;(\mathrm{mod}\;n)$. Hence $X=1$ and (\ref{E:09}) reduces to the equation $D\varphi_2(\lambda)=\varphi_2(\lambda + \gamma h_1)$. In particular $\varphi(\lambda)$ is $\gamma (\omega_i-\omega_j)$-periodic. But by our assumption, $\varphi(\lambda)$ is $(n \omega_i)$-periodic and $\gamma$ is real and irrational. Therefore $\varphi(\lambda)$ is constant and it is a morphism in the category $\CB$.$\square$\\
Now, considering $x$ as a parameter, we obtain:
\begin{cor} The functor $\mathcal{H}^c:\CB^{x} \to \DB^{x}$ is fully faithful. \end{cor}
\paragraph{Remark:} equation (12) shows that $\mathrm{Hom}_{\DB} (I_x^c,I_x^c \otimes V)=\mathrm{Hom}_{\CB}(V^*,I_{x+c}^0)$. Thus the above proposition states that for any finite-dimensional representation $V \in \mathcal{O}b(\CF)$ and for all but finitely many $x\; \mathrm{mod}\; \Z + \tau \Z$, we have $\mathrm{Hom}_{\DB}(V^*,I^0_{x})=\mathrm{Hom}_{\CB}(V^*,\C)$, where the somorphism is induced by the embedding $\C \subset I_x^0$ (constant functions). However, for finitely many values of $x\; \mathrm{mod}\\\; \Z + \tau \Z$, this may not be true: see \cite{KrZ} and \cite{H2} where some finite-dimensional subrepresentations of $I^0_{x}$ are considered.
\section{Full faithfullness of the functor $\mathcal{F}^c_{x}:\CF \to \DB$}
\paragraph{}In this section, we prove the following result:
\begin{prop}The functor $\mathcal{F}^c_{x}: \CF \to \DB$ is fully faithful.
\end{prop}
\textbf{Proof:} we have to show that for any two objects $V,V'$ in $\CF$ there is an isomorphism $\mathcal{F}^c_{x}: \;\mathrm{Hom}_{\CF}(V,V') \to \mathrm{Hom}_{\DB}(\mathcal{F}^c_{x}(V),\mathcal{F}^c_{x}(V'))$. Since $\mathcal{F}^c_{x}$ is trivial at the level of morphisms, this map is injective. Now let $V,W \in {\mathcal{O}}b(\CF)$ and let $\varphi(\lambda) \in \mathrm{Hom}_{\DB}(\mathcal{F}^c_{x}(V),\mathcal{F}^c_{x}(W))$. By definition, $\varphi(\lambda):V \to W$ satisfies the relation
\begin{align*}
\varphi_2(\lambda) S_1(z-x,\lambda-&\gamma h_2)L^V_{12}(z,\lambda)e^{-\gamma D_1} S_1(z-x-c,\lambda)^{-1}\\
&=S_1(z-x,\lambda-\gamma h_2)L^{W}_{12}(z,\lambda)e^{-\gamma D_1} S_1(z-x-c,\lambda)^{-1} \varphi_2(\lambda)
\end{align*}
where $L^V(z,\lambda)$ (resp. $L^W(z,\lambda)$) is the L-operator of $V$ (resp. $W$). This is equivalent to
\begin{equation}\label{E:1000}
\varphi_2(\lambda)S_1(z-x,\lambda-\gamma h_2)L^V_{12}(z,\lambda)=S_1(z-x,\lambda -\gamma h_2) L^W_{12}(z,\lambda)\varphi_2(\lambda-\gamma h_1)
\end{equation}
\hbox to1em{\hfill}Introduce the following notations: write $W=\bigoplus_{\xi} W_\xi$, $V=\bigoplus_\mu V_\mu$, $\varphi(\lambda)=\sum_\nu \varphi_\nu(\lambda)$ for the weight decompositions (so that $\varphi_\nu: V_{\xi} \to W_{\xi +\nu}$). Also let $S(z-x,\lambda)=\sum_{i,j}S^{ij}(z-x,\lambda)E_{ij}$, $L^V_{12}(z,\lambda)=\sum_{i,j}E_{ij}\otimes L_V^{ij}(z,\lambda)$ and use the same notation for $L^W(z,\lambda)$. Applying (\ref{E:1000}) to $v_i \otimes \zeta_\mu$ for some $i$ and $\zeta_\mu \in V_\mu$ yields
\begin{equation}\label{E:101}
\begin{split}
\sum_{j,k,\nu}& S^{kj}(z-x,\lambda-\gamma(\mu+\omega_i-\omega_j))v_k \otimes \varphi_\nu(\lambda)(L_V^{ji}(z,\lambda)\zeta_\mu)\\
&=\sum_{l,k,\sigma} S^{kl}(z-x,\lambda-\gamma(\mu + \omega_i-\omega_l+\sigma)) v_k \otimes L_W^{li}(z,\lambda) \varphi_\sigma(\lambda-\gamma \omega_i)\zeta_\mu
\end{split}
\end{equation}
where we used the weight-zero property of $L^V(z,\lambda)$ and $L^W(z,\lambda)$. Applying $v_k^*$ to (\ref{E:101}) and projecting on the weight space $W_{\mu+\omega_i+\xi}$ gives the relation
\begin{equation}\label{E:102}
\begin{split}
\sum_{\overset{\nu,j}{\nu-\omega_j=\xi}}& S^{kj}(z-x,\lambda-\gamma(\mu+\omega_i-\omega_j))\varphi_\nu(\lambda)(L_V^{ji}(z,\lambda)\zeta_\mu)\\&=\sum_{\overset{\sigma,l}{\sigma-\omega_l=\xi}} S^{kl}(z-x,\lambda-\gamma(\mu+\omega_i-\omega_j+\sigma))L^{li}_W(z,\lambda)(\varphi_\sigma(\lambda-\gamma\omega_i)\zeta_\mu)
\end{split}
\end{equation}
for any $i,k,\xi$ and $\zeta_\mu \in V_\mu$. Now let $A=\{\chi\;|\;\varphi_\chi(\lambda) \neq 0\}$. Fix some $j$ and let $\beta \in A$ be an extremal weight in the direction $-\omega_j$ (i.e $\beta-\omega_j+\omega_k \not\in A$ for $k \neq j$). Then (\ref{E:102}) for $\xi=\beta-\omega_j$ reduces to
\begin{equation}\label{E:103}
\begin{split}
S^{kj}(z-x,\lambda-&\gamma(\mu+\omega_i-\omega_j))\varphi_\beta(\lambda)(L^{ji}_V(z,\lambda)\zeta_\mu)\\&=S^{kj}(z-x,\lambda-\gamma(\mu+\omega_i-\omega_j+\beta)) L^{ji}_W(z,\lambda)\varphi_\beta(\lambda-\gamma \omega_i)\zeta_\mu
\end{split}
\end{equation}
\textbf{Claim:} there exists $i\in \{1,\ldots n\}$, $\mu$ and $\zeta_\mu \in V_\mu$ such that $\varphi_\beta(\lambda)(L^{ji}_V(z,\lambda)\zeta_\mu)\neq 0$ for generic $z$ and $\lambda$.\\
\textbf{Proof:} recall the central element $\mathrm{Qdet}(z,\lambda)\in \mathcal{E}_{\tau,\frac{\gamma}{2}}(\mathfrak{gl}_n)$. By definition, its action on $V$ is invertible. Expanding $\mathrm{Qdet}(z,\lambda)$ along the $j^{th}$-line, we have $\mathrm{Qdet}(z,\lambda)=\sum_i L^{ji}_V(z,\lambda)P_i(z,\lambda)$ for some operators $P_i(z,\lambda) \in \mathrm{End}(V)$. In particular, $\sum_i\mathrm{Im}\;L^{ji}(z,\lambda)=V$, and the claim follows.\\
Thus, the ratio $S^{kj}(z-x,\lambda-\gamma (\mu + \omega_i-\omega_j+\beta))/S^{kj}(z-x,\lambda-\gamma(\mu+\omega_i-\omega_j))$ is independent of $k$. This is possible only if $\beta \in \sum_{r \neq j}\C E_{rr}^*$. Applying this to $j=1,\ldots n$, we see that $A=\{0\}$. Hence $\varphi(\lambda)$ is an $\h$-module map. But then relation (\ref{E:1000}) reduces to $\varphi_2(\lambda)L^V_{12}(z,\lambda)=L_{12}^W(z,\lambda)\varphi_2(\lambda-\gamma h_1)$, and $\varphi(\lambda)$ is an intertwiner in the category $\CF$.$\square$
\begin{cor}The functor $\mathcal{F}^c:\CF^{x} \to \DB^{x}$ is fully faithful. \end{cor}
\section{The image of the vector representation}
\paragraph{}Let us denote $\tilde{V}_F(w)=(\C^n, \chi(w)R^F(w,\lambda))$. It is an object of $\CF$ which equals the tensor product of the vector representation $V_F(w)$ by the one-dimensional representation $(\C,\chi(z))$.
\begin{prop} For any $x,w,\;x+c\not\equiv w\;(\mathrm{mod}\; \Z + \tau\Z)$, we have $\mathcal{F}^c_{x}(V_F(w))\simeq \mathcal{H}^c_{x}(V_B(w))$.
\end{prop}
\textbf{Proof:} by definition, we have 
$$\mathcal{F}^c_{x}(\tilde{V}_F(w))=(\C^n,\chi(z)S_1(z-x,\lambda-\gamma h_2)R^F(z-w,\lambda)e^{-\gamma D_1} \times S_1(z-x-c,\lambda)^{-1}),$$
$$I^c_{x} \otimes V_B(w)=(\C^n,\chi(z)S_1(z-x,\lambda)e^{-\gamma D_1}S_1(z-x-c,\lambda)R^B(z-w))$$
 We claim that the map $\varphi(\lambda)=e^{-\gamma D} (S(w-x-c,\lambda)^{-1})e^{\gamma D} \in \mathrm{End}(\C^n)$ is an intertwiner $\mathcal{H}^c_{x}(V_B(w)) \simeq I^c_{x} \otimes V_B(w) \stackrel{\sim}{\to} \mathcal{F}^c_{x}(\tilde{V}_F(w))$. Indeed, we have
\begin{align*}
S_1&(z-x,\lambda-\gamma h_2)R^F(z-w,\lambda)e^{-\gamma D_1}S_1(z-x-c,\lambda)^{-1}(1 \otimes \varphi(\lambda))\\
&=e^{-\gamma D_2}S_1(z-x,\lambda)e^{\gamma D_2}R^F(z-w,\lambda)e^{-\gamma(D_1+D_2)}S_1(z-x-c,\lambda + \gamma h_2)^{-1}S_2(w-x-c,\lambda)^{-1}e^{\gamma D_2}\\
&=e^{-\gamma D_2}S_1(z-x,\lambda)e^{-\gamma D_1}R^F(z-w,\lambda)S_1(z-x-c,\lambda + \gamma h_2)^{-1}S_2(w-x-c,\lambda)^{-1}e^{\gamma D_2}\\
&=e^{-\gamma D_2}S_1(z-x,\lambda)e^{-\gamma D_1}S_2(w-x-c,\lambda+\gamma h_1)^{-1}S_1(z-x-c,\lambda)^{-1}R^B(z-w)e^{\gamma D_2}\\
&=e^{-\gamma D_2}S_1(z-x,\lambda)S_2(w-x-c,\lambda)e^{-\gamma D_1}S_1(z-x-c,\lambda)^{-1}R^B(z-w)e^{\gamma D_2}\\
&=(1 \otimes \varphi(\lambda)) S_1(z-x,\lambda)e^{-\gamma D_1}S_1(z-x-c,\lambda)^{-1}R^B(z-w)
\end{align*}
where we used Lemma 2 and the zero-weight property of $R^F(u,\lambda)$.$\square$
\begin{lem}Let $V,V' \in {\mathcal{O}}b(\CF)$, $W,W' \in {\mathcal{O}}b(\CB)$ and suppose that $\mathcal{F}^c_{x}(V) \simeq \mathcal{H}^c_{x}(W)$ and $\mathcal{F}^c_{x}(V') \simeq \mathcal{H}^c_{x}(W')$. Then $\mathcal{F}^c_{x}(V \otimes V') \simeq \mathcal{H}^c_{x}(W \otimes W')$.
\end{lem}
\textbf{Proof:} If $\varphi(\lambda):V \to W$ and $\varphi'(\lambda):V' \to W'$ are intertwiners then it is easy to check using the methods above that $\varphi_2'(\lambda-\gamma h_1)\varphi_1(\lambda): V \otimes V' \to W \otimes W'$ is an intertwiner.$\square$
\paragraph{}Applying this to tensor products of the vector representations, we obtain
\begin{cor} For any $x \in \C$ and $w_1,\ldots,w_r \in \C\backslash \{x+c + \Z + \tau\Z\}$, we have
$$\mathcal{F}^c_{x}(\tilde{V}_F(w_1) \otimes \ldots \tilde{V}_F(w_r)) \simeq \mathcal{H}^c_{x} (V_B(w_1) \otimes \ldots V_B(w_r)).$$
\end{cor}
\begin{cor} For any $w_1,\ldots,w_r \in \C$, we have
$$\mathcal{F}^c(\tilde{V}_F(w_1) \otimes \ldots \tilde{V}_F(w_r)) \simeq \mathcal{H}^c (V_B(w_1) \otimes \ldots V_B(w_r)).$$
\end{cor}
\paragraph{}Notice that in this case, we have a canonical intertwiner, given by the formula
$$\varphi_{1\ldots r}(\lambda, w_1, \ldots ,w_r)=\tilde{S}_r^{-1}(w_r-x-c,\lambda-\gamma \sum_{i=1}^{r-1}h_i) \ldots \tilde{S}_{1}^{-1}(w_1-x-c,\lambda),$$
where we set $\tilde{S}(z,\lambda)=e^{-\gamma D}S(z,\lambda)e^{\gamma D}$.
\section{Equivalence of subcategories}
\paragraph{}Let us summarize the results of sections 4-8. By proposition 2, we can identify $\CB^{x}$ with a full subcategory $\mathcal{D}_1^{x}$ of $\DB^{x}$. By proposition 3, we can identify $\CF^{x}$ with a full subcategory $\mathcal{D}_2^{x}$ of $\DB^{x}$. Moreover, $\mathcal{D}_1^{x}$ and $\mathcal{D}_2^{x}$ intersect (at least if we replace $\DB^{x}$ by the equivalent category $\widetilde{\DB^{x}}$ whose objects are isomorphism classes of objects of $\DB^{x}$), and the intersection contains objects of the form $\mathcal{F}^c(\bigotimes_i V_F(w_i))\simeq \mathcal{H}^c(\bigotimes_i V_B(w_i))$, where $i=1, \ldots r$ and $w_i \in \C$. Hence,
\begin{theo} The abelian subcategory $\mathcal{V}_B^{x}$ of $\CB^{x}$ generated by objects 
$\bigotimes_i V_B(w_i)$ for $i=1, \ldots r$, $r \in \N$ and $w_i \in \C$ and the abelian subcategory $\mathcal{V}_F^{x}$ of $\CF^{x}$ generated by objects 
$\bigotimes_i V_F(w_j)$ for $j=1, \ldots s$, $s \in \N$ and $w_j \in \C$ are equivalent.
\end{theo}
\paragraph{}Note that for numerical values of $x$, $\mathcal{F}^c_{x}:\CF \to \DB$ is always fully faithful, and $\mathcal{F}_x^c(\CF)$ a full subcategory of $\DB$, but this is not true of $\mathcal{H}^c_x$, because of the existence of nontrivial finite-dimensional subrepresentations of $I^0_{x}$.
\paragraph{Acknowledgments:} the authors were supported by the NSF grant DMS-9700477. O.S would like to thank Harvard University Mathematics Department for its hospitality without which this work would not have been possible.
\small{ 
}
\noindent Olivier Schiffmann, ENS Paris, 45 rue d'Ulm 75005 Paris, FRANCE \\
\texttt{schiffma@clipper.ens.fr},\\
Pavel Etingof, Harvard Mathematics Dept., Harvard University, Cambridge MA 02138 USA \\
\texttt{etingof@math.harvard.edu} 
\end{document}